

\input amstex.tex
\documentstyle{amsppt}

\define\ZZ{{\Bbb Z}}
\define\NN{{\Bbb N}}
\define\RR{{\Bbb R}}
\define\QQ{{\Bbb Q}}

\define\M{{\Cal M}}

\define\Ha{{\Cal H}}

\define\geg{{\goth g}}
\define\hh{{\goth h}}

\magnification1200

\topmatter

\title
A lecture on Kac--Moody Lie algebras of the arithmetic type
\endtitle

\author
Viacheslav V. Nikulin \footnote{Partially supported by
Grant of Russian Fund of Fundamental Research;
Grant of AMS; and Grant of ISF MI6000.\hfill\hfill}
\endauthor

\address
Steklov Mathematical Institute,
ul. Vavilova 42, Moscow 117966, GSP-1, Russia.
\endaddress

\email
slava\@nikulin.mian.su
\endemail

\abstract
We name an indecomposable symmetrizable
generalized Cartan matrix $A$ and the corresponding
Kac--Moody Lie algebra ${\goth g} ^\prime (A)$
{\it of the arithmetic type} if
for any $\beta \in Q$ with $(\beta | \beta)<0$ there exist
$n(\beta )\in {\Bbb N}$ and an imaginary root $\alpha \in \Delta^{im}$
such that
$n(\beta )\beta \equiv \alpha \mod Ker\ (.|.)$ on $Q$. Here $Q$ is the root
lattice. This generalizes "symmetrizable hyperbolic" type of Kac and Moody.

We show that generalized Cartan matrices of the arithmetic type are divided in
$4$ types: (a) finite, (b) affine, (c) rank two, and (d) arithmetic hyperbolic
type. The last type is very closely related with arithmetic groups generated
by reflections in hyperbolic spaces with the field of definition $\Bbb Q$.

We apply results of the author and \'E.B. Vinberg on arithmetic groups
generated by reflections in hyperbolic spaces to describe generalized
Cartan matrices of the arithmetic hyperbolic type and to show that there
exists a finite set of series of the generalized Cartan matrices of
the arithmetic hyperbolic type.
For the symmetric case all these series are known.
\endabstract

\rightheadtext{Kac--Moody Lie algebras of the arithmetic type}

\leftheadtext{Viacheslav V. Nikulin}

\endtopmatter

\document

\subhead
0. Introduction
\endsubhead

This lecture was given by the author at Johns Hopkins University,
Notre Dame University, Penn State University and Queen's University on
the fall 1994. I am grateful to these Universities for hospitality.

I am  grateful to Professor V. Chari for useful discussions.
I am grateful to Professor \'E.B. Vinberg for his interest to this
subject.

\smallpagebreak

We want to pay attention to one class of Kac--Moody Lie algebras which is
very closely related with arithmetic reflection groups in hyperbolic
spaces.

\subhead
1. Reminding on symmetrizable Kac--Moody Lie algebras
\endsubhead

Here we recall results on symmetrizable Kac--Moody Lie algebras which we need.
One can find them in the book by Victor Kac \cite{Ka1}.

\noindent
(1.1) An $n \times  n$-matrix $A=(a_{ij})$ is called a
{\it generalized Cartan matrix} if

\noindent
(C1)\ \ \  $a_{ii}=2$ for  $i=1,...,n$;

\noindent
(C2) $a_{ij}$ are non-positive integers for $i\not=j$;

\noindent
(C3) $a_{ij}=0$ implies $a_{ji}=0$.

\noindent We denote by $l$ the rank of $A$ and by $k=n-l$ the dimension of
the kernel of $A$.

For simplicity, below we suppose that $A$ is {\it indecomposable} which means
that there does not exist a decomposition $I=\{1,...,n\}=I_1 \cup I_2$
such that both $I_1$ and $I_2$ are non-empty and $a_{ij}=0$ for any
$i \in I_1$ and any $j \in I_2$.

A generalized Cartan matrix $A$ is called {\it symmetrizable} if
there exists an invertible diagonal matrix $D=\text{diag}(\epsilon_1,...,
\epsilon_n)$ and a symmetric matrix $B=(b_{ij})$, such that
$$
A=DB;\ \ \text{or}\ \ (a_{ij})=(\epsilon_ib_{ij})
\tag1.2
$$
One always can suppose that
$$
\epsilon_i \in \QQ,\ \ \epsilon_i>0\ \text{for any $1\le i \le n$},\ \
b_{ij} \in \ZZ\  \text{and}\  B.C.D.(\{ b_{ij}\ |\ 1\le i,j \le n\ \} )=1.
\tag1.3'
$$
By (C1)---(C2), this is equivalent to
$$
b_{ij} \in \ZZ,\ b_{ii}>0,\ b_{ij}\le 0\
\text{for}\ i\not=j\  \text{and}\  B.C.D.(\{ b_{ij}\ |\ 1\le i,j \le n\ \} )=1.
\tag1.3
$$
Then the matrices $D$ and $B$ are defined uniquely. Later we always
suppose that these conditions (1.3) are satisfied.

One formally  defines the {\it root lattice} and the {\it root semigroup}
$$
Q=\bigoplus_{i=1}^{n}{\ZZ \alpha_i},\ \ \ Q_+= \sum_{i=1}^{n} {\ZZ_+\alpha_i}
\subset Q
\tag1.4
$$
where $\ZZ_+$ denotes non-negative integers.
Similarly, one defines the {\it coroot lattice} $Q^\vee$ and the
{\it coroot semigroup} $Q^\vee_+$.
$$
Q^\vee=\bigoplus_{i=1}^{n}{\ZZ \alpha_i^\vee},\ \ \
Q_+^\vee = \sum_{i=1}^{n} {\ZZ_+\alpha_i^\vee} \subset Q^\vee .
\tag1.5
$$
One has a natural pairing
$$
\langle . | . \rangle \ :\ Q^\vee \times Q \to \ZZ,\ \
\langle \alpha_i^\vee , \alpha_j \rangle
=a_{ij}\ \ (i,j =1,...,n)
\tag1.6
$$
and an integral symmetric bilinear form
$$
(. |. )\ :\ Q \times Q \to \ZZ , \ \ (\alpha_i|\alpha_j)=b_{ij}=
a_{ij}/\epsilon_i .
\tag1.7
$$
The symmetric bilinear form $(.|.)$ is called {\it canonical}.
Pairings (1.6) and (1.7) are connected by the formula
$$
a_{ij}=\langle \alpha_i^\vee , \alpha_j \rangle =
2(\alpha_i|\alpha_j)/(\alpha_i|\alpha_i).
\tag1.8
$$
{\it Kac--Moody Lie algebra} $\geg ^\prime (A)$ is the complex
Lie algebra defined by
$3n$ generators $e_1,...,e_n$, $f_1,...,f_n$, $\alpha_1^\vee,...,
\alpha_n^\vee$ and defining relations
$$
\split
&[\alpha_i^\vee, \alpha_j^\vee ]=0 \\
&[e_i,f_j]=\delta_{ij}\alpha_i^\vee,\\
&[\alpha_i^\vee, e_j]=a_{ij}e_j,\\
&[\alpha_i^\vee,f_j]=-a_{ij}f_j,\\
&(ad\ e_i)^{1-a_{ij}}e_j=0,\ \ (ad\ f_i)^{1-a_{ij}}f_j=0\ \text{if}\ i\not=j,
\endsplit
\tag1.9
$$
for any $1\le i,j \le n$, where $\delta_{ij}$ is the Kroneker symbol. The
main property of $\geg^\prime (A)$ is that $\hh^\prime =Q^\vee \otimes
{\Bbb C}$ is the maximal commutative subalgebra of $\geg^\prime (A)$; the
center of $\geg^\prime (A)$ is equal to the
$$
{\goth c}=\{x \in \hh^\prime\ |\ \langle x, Q \rangle =0 \},
\tag1.10
$$
and the Lie algebra $\geg^\prime (A)/{\goth c}$ is simple. The Lie algebra
$\geg ^\prime (A)$ has so called {\it root space decomposition}
$$
\geg ^\prime (A)=
(\bigoplus_{\alpha \in Q_ +}{ \geg_{-\alpha}}) \oplus \hh^\prime \oplus
(\bigoplus_{\alpha \in Q_ +}{ \geg_\alpha})
\tag1.11
$$
where for a non-zero $\alpha \in \pm Q_+$, one set
$$
\geg_{\alpha}=\{ x\in \geg ^\prime (A) \ | \
[h,x]=\langle h, \alpha \rangle x,\ \forall\ h\in \hh ^\prime \}.
$$
It is known that $\dim \geg_{\alpha} < \infty$.
An element $0 \not= \alpha \in \pm Q_+$ is called a
{\it root} if $\dim \geg_\alpha > 0$;
the dimension $\dim \geg_\alpha$ is called the {\it multiplicity} of
the root $\alpha$. The set of roots is denoted by
$\Delta$. It is a disjoint union $\Delta =\Delta_+ \cup -\Delta_+$ where
$\Delta_+ = \Delta \cap Q_+$. A description of
roots and their multiplicities is an important problem for the theory
of Kac--Moody Lie algebras.

Due to Victor Kac, we have the following description of
the set of roots $\Delta$. Evidently,
$\alpha_1,...,\alpha_n$ are roots (they are called {\it simple roots}).
One defines {\it fundamental reflections}
$r_{\alpha_i} \in GL(Q)$
by the formula
$$
r_{\alpha_i}(x)=x-(2(\alpha_i|x)/(\alpha_i|\alpha_i))\alpha_i=
x-\langle \alpha_i^\vee, x\rangle \alpha_i,\ \ x \in Q,
\tag1.12
$$
and {\it Weyl group} $W$ generated by all reflections $r_{\alpha_i}$,
$1\le i \le n$. Evidently, $W$ preserves the canonical symmetric
bilinear form $(.|.)$.
The set of roots $\Delta$ is invariant with respect to
$W$ and is divided on two parts:
{\it real roots} $\Delta^{re}$ and
{\it imaginary roots} $\Delta^{im}$. Here, by definition,
$$
\Delta^{re}=W(\alpha_1)\cup...\cup W(\alpha_n).
\tag1.13
$$
Let
$$
K=\{\alpha \in Q_+-\{0\}\ | \
(\alpha | \alpha_i) \le 0\
\text{for all $\alpha_i,\ 1\le i \le n$, and $supp\ \alpha$ is connected}\}.
\tag1.14
$$
Here, for $\alpha = \sum_{i=1}^{i=n}{k_i\alpha_i} \in Q_+$, the
subset $supp\ \alpha \subset \{\alpha_1,...,\alpha_n\}$
is the set of all $\alpha_i$ such that
$k_i>0$. The $supp\ \alpha$ is called connected if there does not
exist a decomposition $supp\ \alpha =A_1\cup A_2$ such that $(A_1|A_2)=0$,
with non-empty $A_1$ and $A_2$. One has
$$
\Delta^{im}_+=W(K).
\tag1.15
$$
Since $(.|.)$ is $W$-invariant, from these results, it follows that
$$
(\alpha|\alpha)>0,\ \text{if $\alpha \in \Delta^{re}$},
\tag1.16
$$
and
$$
(\alpha|\alpha)\le 0,\ \text{if $\alpha \in \Delta^{im}$}.
\tag1.17
$$
Moreover, if $\alpha \in \Delta^{im}$, then $n\alpha \in \Delta^{im}$ for
any $n \in \NN$.

\subhead
2. Kac--Moody Lie algebras of the arithmetic type
\endsubhead

\definition{Definition 2.1}
A generalized Cartan matrix $A$ and the corresponding
Kac -- Moody Lie algebra $\geg^\prime (A)$ have the {\it arithmetic type}
if $A$ is symmetrizable indecomposable and
for the corresponding canonical
symmetric bilinear form
$(.|.)$ on the root lattice $Q$ one has:
for each $\beta \in Q$ with the property $(\beta | \beta) < 0$ there
exists $n(\beta) \in \NN$ and an imaginary root
$\alpha \in \Delta^{im}$
such that
$$
n(\beta)\beta\equiv \alpha \mod~Ker~(.|.)\ \text{on}\ \ Q.
\tag2.1
$$
Thus, the inequality $(\beta |\beta )< 0$ on the root lattice $Q$
should define imaginary roots up to the kernel of $(.|.)$ on
$Q$ and multiplying by natural numbers.
\enddefinition

To formulate results, let us denote by
$$
S:M \times M \to \ZZ
\tag2.2
$$
the induced by $(.|.)$ canonical
non-degenerate integral symmetric bilinear form on the free
$\ZZ$-module $M=Q /\text{Ker\ }(.|.)$. We denote by $\pi:Q\to M$
the corresponding factorization map. To be shorter, we sometimes
denote $\tilde x= \pi (x)$.
In particular, we denote $\widetilde{W}\subset O(S)$ the image of
$W$ by $\pi$.
Thus, $S$ is the symmetric bilinear form
defined by the integral symmetric matrix $B$ (see (1.2)) modulo
its kernel. We denote by $(t_{+}, t_{-}, t_{0})$ the signature
of a symmetric matrix. Thus, $t_{+}$, $t_{-}$ and $t_{0}$ are equal to
numbers of positive, negative and zero "squares" respectively.

We have the following basic result which is
well-known for the first cases (a), (b) and (c).

\proclaim{Theorem 2.1} A symmetrizable indecomposable generalized
Cartan matrix $A$ and the corresponding Kac--Moody Lie algebra
$\geg ^\prime (A)$ have the arithmetic type if and only if
$A$ has one of the
types (a), (b), (c) or (d) below:

(a) The finite type case: $B>0$ (equivalently, $B$
has the signature $(l,0,0)$).

(b) The affine type case: $B\ge 0$ and $B$ has a $1$-dimensional kernel
(equivalently,
$B$ has the signature $(l,0,1)$).

(c) The rank two hyperbolic case:  $B$ is hyperbolic of the rank $2$
(equivalently, $B$ or $S$ have the signature $(1,1,0)$).

(d) The arithmetic hyperbolic case:
$B$ is hyperbolic of the rank $>2$ (equivalently, $B$ has the
signature $(l-1,1,k)$ where $l\ge 3$, or $S$ has the signature $(l-1,1)$
where $l\ge 3$) and  the index $[O(S):\widetilde{W} ]$ is finite.
\endproclaim

\demo{Proof} Let  $(t_{+}, t_{-}, t_{0})$ be the signature of $B$. If
$t_{-}=0$, we get cases (a) and (b), this is well-known (see \cite{Ka1}).
If $t_{-}=1$ and $t_{+}=1$, we get the case (c). This is also well-known
(see \cite{Ka1}).

Thus, we assume that $t_{-}\ge 1$ and $l=t_{+}+t_{-}\ge 3$.
Let
$$
\RR_+Q_+\subset Q\otimes \RR
$$
be the corresponding real cone in $Q\otimes \RR$ generated by $Q_+$ with the
origin at $0$, and
$$
\widetilde{\RR_+Q_+}=\RR_+\tilde{\alpha_1}+ \cdots +\RR_+\tilde{\alpha_n }
$$
its projection by $\pi$.
We
claim that
$$
\widetilde{\RR_+Q_+} \cap -\widetilde{\RR_+Q_+}=\{0\}.
\tag2.3
$$
Otherwise, there are real $\lambda_i\ge 0$, $1 \le i \le n$, which are not
all equal to $0$ such that
$s=\lambda_1\alpha_1+ \cdots +\lambda_n \alpha_n \in Ker\ (.|.)$.
By \cite{Ka1, Theorem 5.6}, there exists
$\alpha =\sum_{i=1}^{n}k_i\alpha_i$ such that $k_i>0$ and
$(\alpha|\alpha_i)<0$ for any $1\le i \le n$.
Then evidently, $(s|\alpha)<0$. since all $\lambda_i\ge 0$ and not all of them
are $0$. We get a contradiction with $s \in Ker\ (.|.)$.

Let
$$
V(S)=\{x \in M\otimes \RR \ |\ S(x,x)<0 \} .
$$

Let us suppose that $A$ is of the arithmetic type.
Since $M\otimes \QQ$ is everywhere dense in
$M\otimes \RR$ and the set of roots
$\Delta^{im} \subset Q_+\cup -Q_+$, we get that
$$
V(S)\subset \widetilde{\RR_+Q_+} \cup -\widetilde{\RR_+Q_+}
$$
has at least two connected components by (2.3).
It is well known and very easy to see
that $V(S)$ is connected if $t_{-}>1$, and has two connected components
if $t_{-}=1$.  Thus, we have proven that
$t_{-}=1$. Further, we suppose that this is the case.

Thus, further we suppose that $t_{+}\ge 2$ and $t_{-}=1$. Equivalently,
the $S$ is hyperbolic (i.e. of the signature $(t_+,1,k)$)
of the rank $\ge 3$. Then $V(S)$ is a cone which
is  the union of two convex
half-cones: $V(S)=V^+(S) \cup -V^+(S)$ (later we
will choose  the half-cone $V^+(S)$ canonically). We denote by
$$
L(S)=V^+(S)/\RR_+
$$  the corresponding hyperbolic space. Its point is
a ray $\RR_+x,\ x \in V^+(S)$.  Any element $\delta \in M\otimes \RR$
with $S(\delta, \delta)>0$ defines the  half-space
$$
\Ha^+_\delta =\{
\RR_+x \in L(S)\ | \ S(x,\delta) \le 0\}
$$ bounded by the hyperplane
$$
\Ha_\delta =\{ \RR_+x \in L(S)\ |\ S(x,\delta) = 0 \}.
$$
Then the $\delta$ is called the vector which is orthogonal to the half-space
 and the hyperplane. This is defined uniquely up to multiplication
$\lambda\delta$, $\lambda>0$.    Let
$$
\widetilde{\RR_+Q_+}^\ast =\{x
\in M\otimes \RR \ |  \ S(x, \tilde{\alpha_i})\le 0\}
$$
be the dual
cone to  $\widetilde{\RR_+Q_+}$.

\proclaim{Lemma 2.2} We have:
$$
\QQ_+\tilde{K}=(\widetilde{\RR_+Q_+}\cap \widetilde{\RR_+Q_+}^\ast) \cap
M\otimes \QQ.
$$
\endproclaim

\demo{Proof} We first prove that
$$
\pi(\QQ_+K)=
\pi(\QQ_+
\{\alpha \in Q_+-\{0\}\ | \
(\alpha | \alpha_i) \le 0\
\text{for all $\alpha_i,\ 1\le i \le n$}\}).
\tag2.4
$$
Let $0\not= \beta \in Q_+$ and
$(\beta |\alpha_j)\le 0$ for all $1\le j \le n$. If $supp\ \beta$ is
not connected, there exists the canonical decomposition
$\beta=\beta_1+\cdots +\beta_k$ with $\beta_i \in Q_+-\{0\}$,
$supp\ \beta_i \cap supp\ \beta_j=\emptyset$,
$(\beta_i|\beta_j)=0$ if $i\not=j$, and
$k\ge 2$.
Then evidently $\beta_i \in K$ and $(\beta_i|\beta_i)\le 0$.
Since $t_-=1$ and $k\ge 2$, it follows that
$(\beta_i|\beta_i)=0$ and there are $\lambda_i \in \NN$ such
that $\lambda_1\beta_1\equiv \lambda_2\beta_2 \equiv \cdots
\equiv \lambda_k\beta_k \mod Ker\ (.|.)$. Here we also use (2.3).
It follows (2.4). On the other hand, one can easily check that
$$
\pi(\QQ_+
\{\alpha \in Q_+-\{0\}|
(\alpha | \alpha_i) \le 0\
\text{for all $\alpha_i,\ 1\le i \le n$}\})=
(\widetilde{\RR_+Q_+} \cap \widetilde{\RR_+Q_+}^\ast) \cap
M\otimes \QQ.
$$
It follows the proof of Lemma.
\enddemo

Using Lemma 2.2 and (1.15), we get that
$$
\QQ_+\widetilde{\Delta^{im}_+}=
\widetilde{W}(\widetilde{\RR_+Q_+} \cap \widetilde{\RR_+Q_+}^\ast) \cap
M\otimes \QQ.
\tag2.5
$$
Here $\widetilde{\RR_+Q_+} \cap \widetilde{\RR_+Q_+}^\ast$ is
a convex connected cone which is evidently contained in
$\overline{V^+(S)}$ since for
$x\in \widetilde{\RR_+Q_+} \cap \widetilde{\RR_+Q_+}^\ast$ we
have $S(x,x)\le 0$. It follows that we can choose
the half-cone $V^+(S)$ by the condition
$\widetilde{\RR_+Q_+} \cap \widetilde{\RR_+Q_+}^\ast\subset V^+(S)$.
(In fact, when we proved (2.3),
we have mentioned that the cone
$\widetilde{\RR_+Q_+} \cap \widetilde{\RR_+Q_+}^\ast$
is non-degenerate, i.e. it
contains a non-empty open subset of $M\otimes \RR$.)

Since $M\otimes \QQ$ is everywhere dense in $M\otimes \RR$, we
then get that $A$ is of the arithmetic type if and only if
$$
\overline{\widetilde{W}(\widetilde{\RR_+Q_+} \cap \widetilde{\RR_+Q_+}^\ast)}=
\overline{V^+(S)}.
\tag2.6
$$

Let
$$
O_+(S)=\{ \phi \in O(S)\ |\ \phi (V^+(S))=V^+(S) \}
$$
a subgroup of the index two of the automorphism group of $S$. From the
arithmetic of integral symmetric bilinear forms, it is known that
$O_+(S)$ is discrete in the hyperbolic space $L(S)=V^+(S)/\RR_+$ and
has a fundamental domain of a finite volume. The subgroup
$\widetilde{W}\subset O_+(S)$ is generated by reflections in hyperplanes
$\Ha_{\tilde{\alpha_i}}$, $1\le i \le n$, orthogonal to $\tilde{\alpha_i}$.
By (1.3), we have $S({\tilde{\alpha_i}}, {\tilde{\alpha_j}})\le 0$ if
$i\not=j$. From the theory of groups generated by reflections in hyperbolic
spaces (see \cite{V1}, \cite{V2}, for example), it follows that
$$
\M=\bigcap_{i} {\Ha_{\tilde{\alpha_i}}^+}
\tag2.7
$$
is a fundamental polyhedron of $\widetilde{W}$ in $L(S)$. Evidently,
$$
\M=V^+(S)\cap \widetilde{\RR_+Q_+}^\ast /\RR_+.
\tag2.8
$$
By (2.6), we then get that $A$ is of the arithmetic type of and only if
the embedding
$$
(\widetilde{\RR_+Q_+} \cap \widetilde{\RR_+Q_+}^\ast )/\RR_+ \subset
(V^+(S)\cap \widetilde{\RR_+Q_+}^\ast)/\RR_+=\M
\tag2.9
$$
is equality.

Let us suppose that $A$ is of the arithmetic type. Thus, we
have the equality for the embedding (2.9). It follows that $\M$ has
finite volume since it is a convex envelope of a finite set of points
in $\overline{L(S)}$ corresponding to edges of the finite polyhedral cone
$\widetilde{\RR_+Q_+} \cap \widetilde{\RR_+Q_+}^\ast$. It follows that
the index $[O(S):\widetilde{W}]<\infty$
because $O_+(S)$ has a fundamental
domain of a finite volume.

Now suppose that $\M=
V^+(S)\cap \widetilde{\RR_+Q_+}^\ast/\RR_+$ has finite volume.
Since $\widetilde{\RR_+Q_+}^\ast$ is a finite polyhedral cone, it is true
if and only if
$\widetilde{\RR_+Q_+}^\ast \subset \overline{V^+(S)}$. Considering dual
cones, we then get
$V^+(S)^\ast =V^+(S) \subset \widetilde{\RR_+Q_+}$. Thus, finally we
get the sequence of embedded cones:
$$
\widetilde{\RR_+Q_+}^\ast\subset V^+(S) \subset \widetilde{\RR_+Q_+}.
$$
It follows that the embedding (2.9) is equality, and that $A$ is of the
arithmetic type.
\enddemo

\subhead
3. Invariants of a Kac--Moody Lie algebra of the arithmetic hyperbolic type
\endsubhead

Let $A$ be a generalized Cartan matrix (or the corresponding
Kac--Moody Lie algebra $\geg^\prime (A)$) of the arithmetic hyperbolic type.
In fact, above, we have defined for $A$ several invariants which we
describe more precisely below.

The main invariant is the
isomorphism class of the {\it reflective primitive hyperbolic
non-degenerate integral
symmetric bilinear form} of the rank $l \ge 3$ (using (1.3) and (2.2)):
$$
S:M \times M \to \ZZ.
\tag3.1
$$
Here {\it hyperbolic} means that $S$ has signature $(l-1,1)$;
{\it primitive}
means that $S/k$ is not integral for any integral $k>1$; {\it reflective}
means that $[O(S):W(S)]<\infty$.
Here $W(S)$ denote a subgroup of $O_+(S)$ generated
by all reflections. Here we consider elements of $O_+(S)$ as motions of
the hyperbolic space $L(S)$. An automorphism $\phi \in O_+(S)$ is called
{\it reflection} if $\phi$ acts in $L(S)$ as a reflection with respect to
a hyperplane of $L(S)$. One can easily see that every reflection
$\phi\in O(S)$ is
equal to $r_\delta$ for some $\delta\in M$ with the property
$S(\delta, \delta )>0$, where
$$
r_\delta\ :\ x \to x-(2S(x,\delta)/S(\delta,\delta))\delta,\ \ x\in M,
\tag3.2
$$
and $r_\delta \in O(S)$ if and only if
$$
(2S(x,\delta)/S(\delta,\delta))\delta \in M\ for\ any\ x\in M.
\tag3.3
$$
(The automorphism $r_\delta$ of $L(S)$ is the reflection in the
hyperplane $\Ha_\delta$ which is orthogonal to $\delta$.)
The reflection $r_\delta$ will not change if one replaces $\delta$ by
$\lambda\delta$, $\lambda \in \QQ^\ast$. Thus, we can take $\delta$ to be
primitive in $M$. Then (3.3) is equivalent to
$$
S(\delta, \delta)|2S(M, \delta),\ \text{for\ primitive\ }\delta\in M .
\tag3.4
$$

The next invariant is:
$$
\text{A finite index subgroup generated by reflections\ }
\widetilde{W}\subset W(S).
\tag3.5
$$
This subgroup is defined up to automorphisms of $O(S)$.
Let $\M$ be a fundamental polyhedron of $\widetilde{W}$ and
$$
\M=\bigcap_{\delta \in P(\M)_{pr}} {\Ha^+_\delta}
\tag3.6
$$
where $P(\M)_{pr}$ is the set of primitive vectors of $M$ which are
orthogonal to codimension $1$ faces of $\M$. The subgroup $\widetilde{W}$
evidently has the property:
$$
\text{the set $P(\M)_{pr}$ generates $M$.}
\tag3.7
$$
This property (3.7) is important since not every subgroup
of $W(S)$ of
a finite index and generated by reflections has this property.

The elements $\tilde{\alpha_1},...,\tilde{\alpha_n}$ above (see (2.7))
are orthogonal to faces of $\M$. Thus,
$\tilde{\alpha_i}=\lambda_i \delta_i$
where $\delta_i \in P(\M)_{pr}$, $\lambda_i\in \NN$. Evidently, here
$\delta_i \mapsto \tilde{\alpha_i}=\lambda_i\delta_i$ is the isomorphism
between sets $P(\M)_{pr}$ and $\tilde{\alpha_1},...,\tilde{\alpha_n}$.
Moreover, here we get another invariant of $A$ which is a function
$$
\lambda: P(\M)_{pr} \to \NN,\ \ \delta_i\mapsto \lambda_i.
\tag3.8
$$
This function satisfies two important properties:
$$
\{\lambda(\delta )\delta\ |\ \delta \in P(\M)_{pr}\}\ \text{generates\ }\ M
\tag3.9
$$
(in particular, $B.C.D.(\{ \lambda(\delta)\ |\ \delta \in P(\M) \})=1$)
and
$$
S(\lambda(\delta)\delta, \lambda(\delta)\delta)|
2S(\lambda(\delta^\prime)\delta^\prime, \lambda(\delta)\delta)\ \
\text{for any $\delta, \delta^\prime \in P(\M)$.}
\tag3.10
$$
The last property follows from (1.8) and axioms (C1), (C2).

Thus, we correspond to a generalized Cartan matrix $A$ and
the Kac--Moody Lie algebra $\geg^\prime (A)$ of the
arithmetic hyperbolic type the triplet of
invariants:
$$
(S,\ \widetilde{W}\subset W(S),\ \lambda:P(\M)_{pr}\to \NN)
\tag3.11
$$
satisfying the conditions (3.1)---(3.10).
We evidently have

\proclaim{Theorem 3.1} Invariants (3.11) define the generalized
Cartan matrix $A$ and Kac--Moody Lie algebra $\geg^\prime (A)$
of the arithmetic hyperbolic type by the formula
$$
A=(2S(\lambda(\delta^\prime)\delta^\prime, \lambda(\delta )\delta )/
S(\lambda(\delta)\delta, \lambda(\delta)\delta)),\ \
\delta, \delta^\prime \in P(\M)_{pr}.
\tag3.12
$$

\endproclaim

\demo{Proof} It follows at once from our construction and (1.8).
\enddemo

\subhead
4. On the classification of Kac--Moody Lie algebras of the
arithmetic hyperbolic type
\endsubhead

The basic fact here is the following result. Its first part (a)
was proved by the author \cite{N4}, \cite{N5}, and the second part
by \'E.B. Vinberg \cite{V3}.

\proclaim{Theorem 4.1} (V.V. Nikulin, \cite{N4}, \cite{N5}, and
\'E.B. Vinberg, \cite{V3}) (a) For a fixed rank $rk\ S=l\ge 3$,
the set of all isomorphism classes of
reflective primitive hyperbolic integral
symmetric bilinear forms $S$ is finite.

(b) If $S$ is a reflective hyperbolic integral symmetric
bilinear form, then $rk\ S \le 30$.

In particular, the whole set of isomorphism classes of
reflective primitive hyperbolic integral
symmetric bilinear forms $S$ of the rank $rk\ S \ge 3$ is finite.
\endproclaim

\smallpagebreak

Thus, by Theorem 4.1, in principle, it is possible to describe all
reflective primitive hyperbolic integral symmetric bilinear forms $S$ of
the $rk\ S \ge 3$.

Now let us fix one of reflective primitive hyperbolic integral symmetric
bilinear forms $S:M\times M\to \ZZ$.

\smallpagebreak

The next invariant of the triplet
(3.11) is a subgroup of a finite index
$\widetilde{W}\subset W(S)$ generated by reflections. We evidently have

\proclaim{Proposition 4.2} For a fixed index
$N=[W(S):\widetilde{W}]$,  the set of generated by
reflections subgroups $\widetilde{W}\subset W(S)$ is finite.
\endproclaim

\demo{Proof} A fundamental polyhedron $\M$ of $\widetilde{W}$ is a convex
polyhedron which is
a union of $N$ fundamental polyhedra of $W(S)$. It follows that
the number of possibilities for $\M$
is finite up to the action of $\widetilde{W}$.
\enddemo

Now, let us choose a subgroup $\widetilde{W}\subset W(S)$ generated by
reflections and of a finite index.
Let us choose a fundamental
polyhedron $\M$ of $\widetilde{W}$, and let $P(\M)_{pr}$ be
the set of primitive elements of $M$ which are orthogonal to
codimension $1$ faces of $\M$ and directed outside (i.e. $P(\M)_{pr}$ is a
minimal set of primitive elements of $M$ with the property
(3.6)). Let us additionally suppose that
we have the property (3.7) for $P(\M)_{pr}$. Thus, we require that
$$
P(\M)_{pr}\ \text{generates\ }M.
\tag4.1
$$
This gives some additional restrictions on $\widetilde{W}$ and even on the
form $S$ itself because a fundamental polyhedron $\M_{0}$ of the
$W(S)$ should then have the same property:
$$
P(\M_{0})_{pr}\ \text{generates\ }M.
\tag4.2
$$

By (3.3), we have:
$$
S(\delta ,\delta)|2S(\delta^\prime, \delta ),\
\delta, \delta^\prime \in P(\M)_{pr}.
\tag4.3
$$
By (4.3), Theorem 2.1 and our construction, it evidently follows

\proclaim{Theorem 4.3}  A reflective primitive hyperbolic integral
symmetric bilinear form
$$
S:M \times M \to \ZZ
$$
and a subgroup of a finite index generated by reflections
$$
\widetilde{W}\subset W(S)
$$
satisfying the condition (4.1) (and (4.2)) for a fundamental polyhedron $\M$
of $\widetilde{W}$ canonically define a generalized Cartan matrix of the
arithmetic hyperbolic type
$$
A(S, \widetilde{W})=(2S(\delta^\prime, \delta )/S(\delta ,\delta)),\ \
\delta, \delta^\prime \in P(\M)_{pr},
$$
with the first two invariants (3.11) equal to the
$(S, \widetilde{W}\subset W(S))$.

In particular, for $\widetilde{W}=W(S)$, the generalized Cartan matrix
$A(S)=A(S, W(S))$ is defined canonically by the reflective form $S$ itself.
\endproclaim

\smallpagebreak

Now, let us consider functions
$\lambda: P(\M)_{pr} \to \NN $
satisfying the conditions (3.9), (3.10).
Thus, we require:
$$
S(\lambda(\delta)\delta,\lambda(\delta )\delta)|
2S(\lambda(\delta^\prime)\delta^\prime,\lambda(\delta)\delta ),\
\delta, \delta^\prime \in P(\M)_{pr}.
\tag4.4
$$
and
$$
\{\lambda(\delta )\delta\ |\ \delta \in P(\M)_{pr}\}\ \text{generates\ }\ M
\tag4.5
$$
(in particular, $B.C.D.(\{ \lambda(\delta)\ |\ \delta \in P(\M)_{pr}\})=1$).

We have

\proclaim{Proposition 4.4} The set of functions
$\lambda:P(\M)\to \NN$ which satisfy the conditions (4.4) and (4.5) is
finite.
\endproclaim

\demo{Proof} Since the polyhedron $\M$ has finite volume, for any
two elements $\delta, \delta^\prime \in P(\M)_{pr}$, there exists
a sequence
$\delta=\delta_1,...,\delta_k=\delta^\prime \in P(\M)_{pr}$
such that
$$
S(\delta_i,\delta_{i+1})\not=0,\ \  i=1,...,k-1
$$
(this property is well-known, see \cite{V1}).

Using this property and (4.3), one can easily prove that up to
replacing $\lambda$ by $t\lambda$ where $t \in \QQ$,
there exists only a finite set of functions
$\lambda:P(\M)_{pr} \to \NN$ satisfying the condition (4.4).
By the condition (4.5), the set of possible
$t$ is also finite.
\enddemo

Now let us choose one of functions $\lambda:P(\M)\to \NN$ satisfying
the conditions (4.4) and (4.5).

Then by our construction and Theorem 2.1, we have

\proclaim{Theorem 4.5}  A reflective primitive hyperbolic integral
symmetric bilinear form
$$
S:M \times M \to \ZZ,
$$
a subgroup of a finite index generated by reflections
$$
\widetilde{W}\subset W(S)
$$
satisfying the condition (4.1) (and (4.2)) for a fundamental polyhedron $\M$
of $\widetilde{W}$, and
a function
$$
\lambda :P(\M)_{pr} \to \NN
$$
satisfying the conditions (4.4) and (4.5)
 canonically define a generalized Cartan matrix of the
arithmetic hyperbolic type
$$
A(S, \widetilde{W}, \lambda)=
(2S(\lambda(\delta^\prime)\delta^\prime, \lambda(\delta )\delta )/
S(\lambda(\delta)\delta ,\lambda(\delta)\delta)),\ \
\delta, \delta^\prime \in P(\M)_{pr},
$$
with the invariants (3.1) equal to the
$(S, \widetilde{W}, \lambda)$.

In particular, the case $A(S, \widetilde{W})$ of the Theorem
4.3 corresponds to the case $\lambda=1$, and the case
$A(S)$ corresponds to the case $\widetilde{W}=W(S)$ and
$\lambda =1$.
\endproclaim

Thus, we have described all possible invariants $(3.11)$ of
generalized Cartan matrices of the arithmetic hyperbolic type and
by Theorem 3.1 gave
the description of all generalized Cartan matrices of the arithmetic
hyperbolic type.

There exists
a finite set of series of these matrices corresponding to
a finite set of reflective primitive hyperbolic integral symmetric
bilinear forms of Theorem 4.1. These algebras almost canonically
correspond to these forms since they are constructed by a subgroup
$\widetilde{W}\subset O(S)$ of a finite index (with a finite
set of additional data --- $\lambda$).
Generalized Cartan matrices $A(S)$ and
$A(S, \widetilde{W})$ canonically
correspond to $S$ and a choice of the subgroup $\widetilde{W}$. Thus,
in principle, all information about these Kac---Moody Lie algebras one
can get from the arithmetic of the reflective primitive hyperbolic
integral symmetric bilinear forms $S$ which are described by
Theorem 4.1.

In \cite{Ka1}, there was considered a very particular case of
generalized Cartan matrices $A$ and corresponding Kac--Moody Lie algebras
$\geg ^\prime (A)$ of the arithmetic hyperbolic type. They are called
{\it hyperbolic}. In our notation, symmetrizable hyperbolic case
is exactly  the case when the fundamental polyhedron $\M$  of
$\widetilde{W}$ is a simplex. There exists only a finite
list of these $A$. These $A$ are characterized by the property:
$0\not= \delta \in Q$
is an imaginary root if and only if $(\delta |\delta)\le 0$.

\smallpagebreak

Unfortunately, the complete list of the reflective forms
$S$ of Theorem 4.1 is not known yet.

\subhead
5. Symmetric case
\endsubhead

Let us consider {\it symmetric} generalized Cartan matrices $A$ of
the arithmetic hyperbolic type. Then we put $B=A$, and the subgroup
$\widetilde{W}\subset O(S)$ is a subgroup of the group
$W^{(2)}(S)$ generated by reflections in vectors $\delta\in M$ such
that $S(\delta,\delta)=2$. Thus, the
hyperbolic integral symmetric bilinear form
$S$ should be {\it 2-reflective}, which
means that $[O(S):W^{(2)}(S)]<\infty$.
All these $2$-reflective forms $S$ and fundamental polyhedra
$\M_0$ for $W^{(2)}(S)$ are found
(see \cite{N1}, \cite{N2}, \cite{N3}
and \cite{N6}). For these forms, the maximum $rk\ S=19$.
Thus, one has a description of all series
of {\it symmetric} generalized Cartan matrices
$A$ of the arithmetic hyperbolic type.

\end